\documentclass[11pt]{amsart}

\usepackage[square,compress,comma, numbers,sort]{natbib}
\usepackage[colorlinks=true, citecolor=red, linkcolor=blue]{hyperref}
\usepackage{amsfonts,mathtools}
\usepackage{cleveref}

\allowdisplaybreaks[4]

\usepackage{amsmath}
\usepackage{bbm}
\usepackage{amssymb}
\usepackage{mathtools}

\usepackage{color}

\definecolor{c20}{rgb}{0.,0.7,0.}
\definecolor{c30}{rgb}{0.,0.,1.}
\definecolor{c40}{rgb}{1,0.1,0.7}
\definecolor{c50}{rgb}{1,0,0}
\definecolor{c60}{rgb}{1,0.9,0.1}

\def\x{\vk{x}}

\usepackage{dsfont}

\newcommand{\ve}{\varepsilon}

\newcommand{\abs}[1]{\left\lvert #1 \right\rvert}

\newcommand{\pk}[1]{\mathds{P} \left\{ #1 \right \} }

\newcommand{\R}{\mathbb{R}}

\newcommand{\inr}{\in \R}

\newcommand{\ldot}{,\ldots,}

\newcommand{\BQN}{\begin{eqnarray}}
	\newcommand{\EQN}{\end{eqnarray}}
\newcommand{\BQNY}{\begin{eqnarray*}}
	\newcommand{\EQNY}{\end{eqnarray*}}

\newcommand{\BS}{\begin{sat}}
	\newcommand{\ES}{\end{sat}}
\newcommand{\BT}{\begin{theo}}
	\newcommand{\ET}{\end{theo}}
\newcommand{\BK}{\begin{korr}}
	\newcommand{\EK}{\end{korr}}

\newcommand{\BD}{\begin{de}}
	\newcommand{\ED}{\end{de}}

\newcommand{\BIT}{\begin{itemize}}
	\newcommand{\EIT}{\end{itemize}}
\newcommand{\BDI}{\begin{description}}
	\newcommand{\EDI}{\end{description}}

\newcommand{\BRM}{\begin{remarks}}
	\newcommand{\ERM}{\end{remarks}}

\newcommand{\BEL}{\begin{lem}}
	\newcommand{\EEL}{\end{lem}}

\newtheorem{theo}{Theorem}[section]
\newtheorem{sat}[theo]{Proposition}
\newtheorem{de}[theo]{Definition}
\newtheorem{lem}[theo]{Lemma}

\newtheorem{example}[theo]{Example}
\newtheorem{korr}[theo]{Corollary}
\newtheorem{remark}[theo]{Remark}
\newtheorem{remarks}[theo]{Remarks}
\newtheorem{prop}[theo]{Proposition}

\newcommand{\COM}[1]{}

\def\ve{\varepsilon}

\def\td{\text{\rm d}}
\def\IF{\infty}

\newcommand{\kb}[1]{\boldsymbol{#1}}
\newcommand{\vk}[1]{\kb{#1}}

\def\ve{\varepsilon}

\def\IF{\infty}

\begin{document}
	
	\title{Multidimensional Brownian Risk Models with Random Trend }
	
	\author{Goran Popivoda}
	\address{University of Montenegro, Faculty of Science and Mathematics \\ Podgorica, Montenegro}
		\email{goranp@ucg.ac.me}

	\author{Timofei Shashkov}
	\address{Department of Actuarial Science, Faculty of Business and Economics\\
		University of Lausanne\\
		UNIL-Dorigny, 1015 Lausanne, Switzerland
	}
	\email{timofei98shashkov@gmail.com}
	
	\bigskip
	
	\date{\today}
	\maketitle

	{\bf Abstract:} Let \(\vk B(t)=(B_1(t), \dots,B_d(t))^\top\), \(t\in[0,T]\), \(d\geq 2\) be a \(d\)-dimensional Brownian motion with independent components and let \(\vk \eta=(\eta_1,\dots,\eta_d)^\top\) be a random vector independent of \(\vk B\) such that
	\[
	\pk{\vk K_{1}\leq\vk\eta\leq\vk K_{2}}
	=\pk{K_{11}\leq\eta_1\leq K_{21},\dots,K_{1d}\leq\eta_d\leq K_{2d}}=1,
	\]
	where \(\vk K_1=(K_{11},\dots,K_{1d})^\top\) and \(\vk K_2=(K_{21},\dots,K_{2d})^\top\) are fixed \(d\)-dimensional vectors.
	
	The goal of this paper is to derive asymptotics of
	\[
	\pk{\exists_{t\in[0,T]}: X_1(t)>a_1u,\dots,X_d(t)>a_du}, \ \ \vk X(t)=\left(X_1(t),\dots,X_d(t)\right)^\top
	=A\vk B(t)-\vk\eta t
	\]
	as \(u\to\IF\) under certain restrictions on the random vector \(\vk\eta\) and constants \(a_1,\dots, a_d\).
	
	{\bf Key Words:} multivariate Brownian risk model; Brownian motion; simultaneous ruin probability;
	simultaneous ruin time; ruin time approximation, random trend.
	
	{\bf AMS Classification:} Primary 60G15; secondary 60G70

	\section{Introduction}
    For a given insurance portfolio, let us consider the risk model defined by
	\[
	R(t)= u+ ct - \sigma Z(t), \ \ \ \ \ t\ge 0,
	\]
	where \(Z(t)\) is a centered Gaussian risk process with a.s. continuous sample paths, \(u>0\) is the initial capital, \(c>0\) is the premium rate and \(\sigma>0\) is a constant. 
	
	The ruin probability over the time horizon \([0,T]\), \(T>0\) is defined by
	\begin{align*}
		\widetilde{\psi}_T(u)\coloneqq\pk{ \inf_{t\in [0,T]} R(t) < 0}= \pk{\sup_{t\in [0,T]} (\sigma Z(t)- c t)> u}
	\end{align*}
	and in view of,  e.g., \cite{DeM15, Boukai}, for \(Z\) being a standard Brownian motion,
	for any \(u\geq 0\)
	\begin{align*} 
		\widetilde{ \psi}_T(u) =\Phi\left(-\frac u{\sigma \sqrt{T}} -\frac{c\sqrt{T}}{\sigma}\right)+
		e^{-2cu/\sigma^2}\Phi\left(- \frac u{\sigma \sqrt{T}} +\frac{c\sqrt{T}}{\sigma}\right),
	\end{align*}
    where \(\Phi\) is the distribution function of a \(N(0,1)\) random variable.

    The asymptotic behavior of the classical ruin probability is a commonly investigated topic (see, for example, \cite{Pit96, Piter2, MR3776537,Jasnov,Filomat}). The primary obstacle in analyzing of extreme values of vector-valued Gaussian processes is that the Slepian inequality, a key tool in studying the excursion probabilities of Gaussian processes, cannot be applied to general vector-valued Gaussian processes. Therefore, there is a lack of an universal approach for investigating the excursion probabilities in the vector-valued case.	

    The two-dimensional Brownian motion risk model, a natural limiting model of several general bivariate insurance risk models, was recently studied in \cite{Del2020}. In \cite{mi:18}, the two-dimensional Brownian risk model's simultaneous ruin probability and simultaneous ruin time approximations are obtained as the initial capital increases to infinity, while the paper \cite{Simruin} considers the simultaneous ruin probability of a $d$-dimensional Gaussian process that consists of mutually independent centered Gaussian processes with stationary increments.

    The paper \cite{Paris} derives asymptotic approximations for the Parisian and cumulative Parisian ruin probabilities, as well as the simultaneous ruin time, in the two-dimensional Brownian risk model as the initial capital approaches infinity. 

    This subject's actuality has been confirmed by papers such as \cite{Kor2020,Kor2021, KRYSTECKI2022, KRYSTECKI2023, Filomat}. Some of these papers have examined non-simultaneous multidimensional models.

	Let us now consider  \(d\geq 2\) portfolios with the risk vector-process \(\vk R(t)\) defined by
	\begin{align*}
		\vk R(t)=(R_1(t),\dots,R_d(t))^\top=\vk u+ \vk c t- A\vk B(t), \quad t\ge 0
	\end{align*}
	such that \(R_i(t)\) is the risk process that corresponds to the \(i\)-th portfolio, \(i=1,\dots,d\). Here \(A\) is a non-degenerate  \(d\times d\) matrix, $\vk B(t)=(B_1(t),\dots,B_d(t))^\top$ is a $d$-dimensional Brownian motion with independent components, $\vk u=(u_1,\dots,u_d)^\top$ is a vector of initial capitals and $\vk c=(c_1,\dots,c_d)^\top$ is a vector of premium rates. Even though it is natural to assume in applications that the vector $\vk c$ has only nonnegative components, we do not impose such conditions here. 
	
	In this paper, we study the event of simultaneous bankruptcy of all \(d\) insurance companies. To be more specific, 
	denote
	\[
	\psi_T(\vk u):=\pk{\exists t\in[0,T]: \vk W(t) - \vk ct > \vk u},
	\]
	where \(\vk W(t)=A\vk B(t)\). In \cite[Thm \(2.3\)]{Pdmc}, it was shown that
	\[
	\psi_T(\vk au)\sim \mathcal{C}_{\vk a}\pk{\vk W(T)>\vk au+\vk cT}
	\]
    as \(u\to\IF\), where \(\vk a\) has at least one positive component and \(\mathcal{C}_{\vk a}\in (0,\IF)\) is a constant (see \cite[Thm \(2.3\)]{Pdmc} or \Cref{Th1} of the present paper). Hereafter by $\sim$ we denote asymptotic equivalence as \(u\to\IF\).
	
	So far we have considered \(A\) and \(\vk c\) to be fixed, but the natural question is whether one can obtain asymptotics of a similar form, when \(A\) and \(\vk c\) are assumed to be random. In general this is a challenging problem. However, under some natural assumptions, this problem can be solved, for example, if we assume $\vk c$ to be random (and we below denote it by \(\vk\eta\)), but \(A\) to be fixed. We also assume that \(\vk\eta\) is independent of \(\vk W(t)\). In fact, it turns out that
	\[
	\psi_T(\vk au)\sim \mathcal{C}_{\vk a}\pk{\vk W(T)>\vk au+\vk\eta T},\ \ \text{ as } u\to\infty,
	\]
	where $\mathcal{C}_{\vk a}$ is the same as for the case of a fixed trend and where we assumed that $\vk\eta$ is almost surely a uniformly bounded random vector, i.e., there exist two $d$-dimensional constants $\vk K_1=(K_{11},\dots,K_{1d})^\top$ and $\vk K_2=(K_{21},\dots,K_{2d})^\top$ such that
	\BQN
	\pk{\vk K_1\leq\vk\eta\leq\vk K_2}=\pk{K_{11}\leq\eta_{1}\leq K_{21},\dots,K_{1d}\leq\eta_{d}\leq K_{2d}}=1.
	\EQN
    We organize the paper in the following way: in Section \(2\) we present and discuss the main result (\Cref{Th1}) followed by illustrating examples. In Section 3 we first provide a brief outline of the proof of \Cref{Th1} and then discuss all auxiliary statements in detail.
    
    \section{Main result}
    To formulate our main result, we need to introduce the quadratic programming problem.
    \subsection{Quadratic programming problem} 
    Before formulating the quadratic programming problem, we introduce some auxiliary notation.
    
    Unless otherwise stated, we assume that the bold designations are vector-valued and that the nonbold designations are scalars.
    
    For some set \(F=\{f_1,f_2,\dots,f_l\}\subset\{1,\dots,k\}\) and some \(k\)-dimensional vector \(\vk x\), we denote \(\vk x_F:=(x_{f_1},x_{f_2},\dots,x_{f_l})\), an \(l\)-dimensional vector that is obtained from the vector \(\vk x\) by taking only coordinates with indices from the set \(F\). 
    
    Similarly,  we define \(\Sigma_{F_1F_2}\), where \(\Sigma\) is a \(k\times k\) matrix and \(F_1,F_2\) are subsets of \(\{1,\dots,k\}\) with \(l_1\) and \(l_2\) elements, respectively: the \(l_1\times l_2\) matrix \(\Sigma_{F_1F_2}\) is formed from \(\Sigma\) by taking only coefficients \(\Sigma_{f_1,f_2}\), where \(f_1\in F_1\) and \(f_2\in F_2\).
    
    For example, if 
    \(\Sigma=\begin{pmatrix}a_{11} & a_{12}
    	     \\a_{21} & a_{22}\end{pmatrix}
    \), 
    \(F_1=\{1,2\}\) and \(F_2=\{2\}\), then 
    \(\Sigma_{F_1F_2}=\begin{pmatrix}a_{12} 
    	              \\a_{22}\end{pmatrix}
    \).
    
    The next proposition, which is stated, for example, in \cite[Lem \(2.1\)]{Quad}, discusses the quadratic programming problem.
    \begin{prop}\label{quad}[Quadratic programming problem]
    	Let $\Sigma$ be a $d\times d$ positive definite matrix and let $\vk a \inr^d \setminus (-\IF, 0]^d $.
    	The quadratic programming problem
    	$\Pi_\Sigma(\vk a)$: "minimize $\vk x^T \Sigma^{-1} \vk x$ under constrain $\vk x \ge \vk a$" has a unique solution $\widetilde{\vk a}$, and there exists a unique nonempty
    	index set $I\subseteq \{1\ldot d\}$ with $m\le d$ elements such that
    	$$ \label{eq:IJ}
    	\widetilde{\vk a}_{I}=\vk a_{I} ,  \quad (\Sigma_{II})^{-1} \vk a_{I}>\vk{0}_I,$$
    	and if $I^c :=\{ 1 \ldot d\}\setminus I \not=\emptyset$, then 
    	$$\label{eq:IJ2} \widetilde{\vk a}_{I^c} = \Sigma_{{I^c}I} (\Sigma_{II})^{-1} \vk a_{I}\ge \vk a_{I^c}.$$
    	Next,
    	$$\min_{\vk x \geq \vk a}\vk x^T \Sigma^{-1} \vk x= \widetilde{\vk a}^T \Sigma^{-1}\widetilde{\vk a}=\widetilde{\vk a}_I^T (\Sigma_{II})^{-1}\widetilde{\vk a}_I>0$$
    	and for any $\vk x\in\R^d$
    	$$\label{eq:new}
    	\vk x^T \Sigma^{-1} \widetilde{\vk a}= \vk x_F^T (\Sigma_{FF})^{-1} \vk a_F,$$
    	with \(F\) being an arbitrary subset of $\{1\ldot d \}$ such that $I\subset F$. If $\vk a= (a \ldot a)^\top , a\in (0,\IF)$,
    	then $ 2 \le \abs{I} \le d$, and ${I^c}$ is empty if $\Sigma^{-1} \vk a > \vk{0}$.\\
    	Conversely, if for some nonempty index set  
    	$I \subset \{ 1 \ldot d \}$ we have 
    	$$(\Sigma_{II})^{-1}\vk a_I> \vk 0_I,  \quad  \Sigma_{{I^c}I} (\Sigma_{II})^{-1} \vk a_I \ge \vk a_{I^c},$$
    	then $\widetilde{\vk a}$ with $\widetilde{\vk a}_{I^c}= \Sigma_{{I^c}I} (\Sigma_{II})^{-1} \vk a_I, \widetilde{\vk a}_I= \vk a_I$ is the solution of $\Pi_{\Sigma}(\vk a)$.
    \end{prop}
    \begin{remark}\label{I_c}
    	From \Cref{quad}, it follows that \(\left(\Sigma^{-1}\widetilde{\vk a}\right)_{I^c}=\vk 0_{I^c}\).
    	
    	Indeed, for all \(\vk\x\in\R^d\) by \Cref{quad}, we have
    	\[
    	\vk x_{I}^T(\Sigma^{-1}\widetilde{\vk a})_I+\vk x_{I^c}^T(\Sigma^{-1}\widetilde{\vk a})_{I^c}
    	=\vk x^T\Sigma^{-1}\widetilde{\vk a}=\vk x_{I}^T(\Sigma_{II}^{-1})\widetilde{\vk a}_I.
    	\]
    	Therefore
    	\[
    	\vk x_{I^c}^T(\Sigma^{-1}\widetilde{\vk a})_{I^c}
    	=\vk x_{I}^T(\Sigma_{II}^{-1})\widetilde{\vk a}_I-\vk x_{I}^T(\Sigma^{-1}\widetilde{\vk a})_I.
    	\]
    	Since the right-hand side in the last equality is independent of \(\vk x_{I^c}\), we can conclude that \(\left(\Sigma^{-1}\widetilde{\vk a}\right)_{I^c}=\vk 0_{I^c}\).
    \end{remark}
    \subsection{Main result}
    Here and below all operations and relations for vectors such as $"+", "/", "<",$ and multiplication by scalar constants, are assumed to be componentwise. For example, \(\vk x/\vk y:=\left(x_1/y_1,\dots,x_k/y_k\right)^\top\) for $k$-dimensional vectors $\vk x$ and $\vk y$.
    
    Let $\vk W(t)=A\vk B(t)$, where $A$ is a fixed non-degenerate $d\times d$ real matrix and $\vk B$ is a $d$-dimensional Brownian motion with independent components, $d\geq 2$. Let also \(\Sigma\) be the covariance matrix of \(\vk W\).
    
    Denote by $\langle\vk x,\vk y\rangle$ for $k$-dimensional vectors \(\vk x,\vk y\), the standard scalar product on the Euclidean space $\R^k$, i.e., 
    $\langle\vk x,\vk y\rangle=\sum_{i=1}^{k}x_iy_i.$
    
    Using the self-similarity of Brownian motion without loss of generality, we assume in the sequel that \(T=1\).
    Denote then
    $$\psi_1(\vk a u)=\pk{\exists_{t\in[0,1]}:\vk W(t)-\vk\eta t>\vk au}.$$
    Using notation from \Cref{quad}, we formulate the main result.
    \begin{theo}\label{Th1} Under the above stated conditions, $\vk a\in\R^d\setminus(-\IF, 0]^d$
    	$$\psi_{1}(\vk au)\sim \left(\prod_{i\in I}\lambda_i\right)I_a\pk{\vk W(1)>\vk au+\vk\eta},\text{ as } u\to\IF,$$
    	where $\vk\lambda=\Sigma^{-1}\widetilde{\vk a}$ and 
    	\[
    	I_a=\int_{\R^{|I|}}\pk{\exists_{t\geq 0}:\vk W_I(t)-\vk a_It>\vk x_I}e^{\langle\vk\lambda_I,\vk x_I\rangle}d\vk x_I\in(0,\IF).
    	\]
    \end{theo}
    
    We continue with illustrating examples and discussions. 
    
    \begin{example}\label{d=2}
    Consider the case where $d=2$ and $A=\begin{pmatrix}1 & 0\\\rho & \rho^*\end{pmatrix}$, with $\vk a=(1,a)^\top$, where $a\leq 1$. Then the covariance matrix of the vector \(\vk W(1)=A\vk B(1)=(B_1(1), \rho B_1(1)+\rho^*B_2(1))^\top\)
    is equal to $\begin{pmatrix}1 & \rho\\\rho & 1\end{pmatrix}.$ If $a>\rho$, then $(1,a)$ is an optimal solution to the quadratic programming problem and $I=\{1,2\}$.
    
    If $a\leq\rho$, then $(1,\rho)$ is an optimal solution with $I=\{1\}$.
    
    Therefore, using the result of \Cref{Th1}, we obtain the following result.
    \begin{korr}
    	\label{Th2} If $a>\rho$, then
    	$$\psi_{1}(u,au)\sim I_{a}\lambda_1\lambda_2\pk{\begin{array}{ccc}
    			W_1(1) - \eta_1 > u \\
    			W_2(1)- \eta_2  > au
    	\end{array}},\quad u\to\IF,$$
    	where $$I_{a}=\int_{\R^2}
    	\pk{\exists_{t\geq 0}:\begin{array}{ccc}
    			W_1(t) - t > x_1 \\
    			W_2(t)- at  > x_2
    	\end{array}} e^{ \lambda_1  x_1 + \lambda_2x_2 }\, \td\vk x,$$
    	and \(\lambda_1=(1-a\rho)/(1-\rho^2)\), \(\lambda_2=(a-\rho)/(1-\rho^2)\).
    	
    	If $a\leq\rho$, then \(I_a=2, \lambda_1=1\) and
    	$$\psi_{1}(u,au)\sim 2\pk{\begin{array}{ccc}
    			W_1(1) - \eta_1 > u \\
    			W_2(1)- \eta_2  > au
    	\end{array}},\quad u\to\IF.$$
    \end{korr}
    \end{example}
    Discuss the model, which generalizes the Example \ref{d=2}. 
    \begin{example}\label{Eq-cor}\((\)Equi correlated risk model\()\).
    	Let \(A\) be such that \(\Sigma=AA^T\) with \(\Sigma_{ii}=1\) for all \(i=1,\dots,d\) and \(\Sigma_{ij}=\rho\in (-1/(d-1), 1)\), if \(i\neq j\), \(i,j\in\{1,\dots, d\}\).
    	
    	As it was shown in \cite{Simruin}
    	\[
    	\Sigma^{-1}=\left[I_d-\vk 1\vk 1^T\frac{\rho}{1+\rho(d-1)}\right]\frac{1}{1-\rho},
    	\]
    	where hereafter \(I_d\) is a \(d\)-dimensional identity matrix and \(\vk 1=(1,\dots,1)^T\in\R^d\).
    	
    	Moreover, according to \cite{Simruin} if \(a_1\geq\dots\geq a_d\), \(a_1>0\), then \(|I|=d\) if and only if
    	\BQN\label{iff}
    	a_d>\frac{\sum_{i=1}^da_i}{1+\rho(d-1)}.
    	\EQN
    	Additionally
    	\[
    	\vk\lambda_{I}=\frac{1}{1-\rho}\left[\vk a_I-\rho\frac{\sum_{i\in I}a_i}{1+\rho(|I|-1)}\vk 1_{I}\right].
    	\]
    	Take \(a_1=\dots=a_d=a>0\) and \(\vk \eta=(\eta,\dots,\eta)^T\), where \(\eta\) is some bounded random variable. Then, according to (\ref{iff}), \(I=\{1,\dots, d\}\), and by \Cref{Th1}, we obtain that
    	\[
        \psi_1(\vk au)\sim  \left(\prod_{i=1}^d\lambda_i\right)I_a\pk{\forall i\in\{1,\dots, d\}:W_i(1)>au+\eta},\text{ as } u\to\IF,
    	\]
    	where
    	\[
    	I_a=\int_{\R^d}\pk{\exists_{t\geq 0}:\vk W(t)-\vk at>\vk x}e^{\langle\vk\lambda,\vk x\rangle}d\vk x\in(0,\IF).
    	\]
    \end{example}
    In the following two examples, we show how the asymptotics can be simplified for specific distributions of \(\vk \eta\).
    \begin{example}\label{Bern}
    Let $\vk\eta=(\eta_1,\dots,\eta_d)^{T}$, where $\eta_k$ is a Bernoulli random variable with parameter $p_k\in [0,1]$ for all $k\in\{1,\dots,d\}$, and $\eta_1,\dots,\eta_d$ are jointly independent.
    Then it turns out that 
    $$\pk{\vk W(1)>\vk au+\vk\eta}\sim\prod_{k=1}^{d}(1-p_k)\pk{\vk W(1)>\vk au}$$
    as $u\to\IF$. Indeed, since $\vk\eta$ is independent of $\vk W(1)$, by the law of total probability,
    \begin{align*}
    	\pk{\vk W(1)>\vk au+\vk\eta}
    	&=\sum_{\vk k\in\{0,1\}^d}\pk{\vk W(1)>\vk au+\vk\eta|\vk\eta=\vk k}\pk{\vk\eta=\vk k}
    	\\&=\prod_{k=1}^{d}(1-p_k)\pk{\vk W(1)>\vk au}+\sum_{\vk k\in\{0,1\}^d:\vk k\neq\vk 0}\pk{\vk W(1)>\vk au+\vk k}\pk{\vk\eta=\vk k}&
    \end{align*}
    Observe that
    $$\sum_{\vk k\in\{0,1\}^d:\vk k\neq\vk 0}\pk{\vk W(1)>\vk au+\vk k}\pk{\vk\eta=\vk k}=o\left(\prod_{k=1}^{d}(1-p_k)\pk{\vk W(1)>\vk au}\right),$$
    as $u\to\IF$.
    Indeed, there exists a uniformly in \(\vk k\) bounded from above and below constant $D_2(\vk k)$ such that
    $$\frac{1}{u^{|I|}}\phi_{\Sigma}(\widetilde{\vk a}u+\vk k)\sim D_2(\vk k)\pk{\vk W(1)>\vk au+\vk k},$$
    as $u\to\IF$ uniformly in $\vk k$ (see \Cref{unsim} below). Hence, there exists a constant $r_0>0$ that is independent of \(\vk k\) such that
    \begin{align*}
    	\sum_{\vk k\in\{0,1\}^d:\vk k\neq\vk 0}\pk{\vk W(1)>\vk au+\vk k}\pk{\vk\eta=\vk k}
    	&\leq r_0\sum_{\vk k\in\{0,1\}^d:\vk k\neq\vk 0}\frac{1}{u^{|I|}}\phi_{\Sigma}(\widetilde{\vk a}u+\vk k)
    	\\&=r_0\sum_{\vk k\in\{0,1\}^d:\vk k\neq\vk 0}\frac{1}{u^{|I|}}\phi_{\Sigma}(\widetilde{\vk a}u)e^{-(\widetilde{\vk a}u)^{T}\Sigma^{-1}\vk k-\frac{\vk k^T\Sigma^{-1}\vk k}{2}}
    	\\&\leq r_0\frac{1}{u^{|I|}}\phi_{\Sigma}(\widetilde{\vk a}u)\sum_{\vk k\in\{0,1\}^d:\vk k\neq\vk 0}e^{-(\widetilde{\vk a}u)^{T}\Sigma^{-1}\vk k},
    \end{align*}
    where we used the fact that $\Sigma^{-1}$ is a positive definite matrix.
    By \Cref{I_c} 
    \[
    (\widetilde{\vk a}u)^{T}\Sigma^{-1}\vk k=\vk\lambda_I^{T}\vk k_Iu, \ \text{ where } \ \vk\lambda_I=\left(\Sigma^{-1}\widetilde{\vk a}\right)_I>\vk 0_I.
    \]
    Hence, since $\vk k_I\geq \vk 0_I$, by \Cref{unsim}, we conclude that
    $$ r_0\frac{1}{u^{|I|}}\phi_{\Sigma}(\widetilde{\vk a}u)\sum_{\vk k\in\{0,1\}^d:\vk k\neq\vk 0}e^{-(\widetilde{\vk a}u)^{T}\Sigma^{-1}\vk k}\leq 2^dr_2\frac{1}{u^{|I|}}\phi_{\Sigma}(\widetilde{\vk a}u)e^{-r_1u}=o\left(\pk{\vk W(1)>\vk au}\right).$$
    where 
    $$r_1=\min_{\vk k\in\{0,1\}^d:\vk k\neq\vk 0}\left(\vk \lambda_I^{T}\vk k_I\right)>0$$
    and \(r_2\) is a positive constant. Therefore, it follows that 
    $$\pk{\vk W(1)>\vk au+\vk\eta}\sim\prod_{k=1}^{d}(1-p_k)\pk{\vk W(1)>\vk au},$$
    as $u\to\IF$. According to \Cref{Th1}, we obtain that
    $$\psi_{1}(\vk au)\sim \prod_{i=k}^{d}(1-p_k)\left(\prod_{i\in I}\lambda_i\right)I_a\pk{\vk W(1)>\vk au}$$
    as \(u\to\IF\).
    \end{example}
    \begin{example}
    	Let now \(\vk\eta\sim Uni([0,1]^d)\), i.e. uniformly distributed on \([0,1]^d\) vector-valued random variable. Moreover, to make calculations simpler, assume that \(I=\{1,\dots,d\}\). Observe first that by the law of total probability, we have that
    	\[
    	\pk{\vk W(1)>\vk au+\vk\eta}=\int_{[0,1]^d}\pk{\vk W(1)>\vk au+\vk c}d\vk c.
    	\]
    	Hence, by \Cref{unsim}, we obtain that
    	\begin{align*}
    		\pk{\vk W(1)>\vk au+\vk\eta}
    		&\sim\frac{1}{\left(\prod_{i=1}^{d}\lambda_i\right)u^{d}}\int_{[0,1]^d}\phi_{\Sigma}(\widetilde{\vk a}u+\vk c)d\vk c
    		=\frac{1}{\left(\prod_{i=1}^{d}\lambda_i\right)u^{d}}\int_{[0,1]^d}e^{-\frac{(\widetilde{\vk a}u+\vk c)^T\Sigma^{-1}(\widetilde{\vk a}u+\vk c)}{2}}d\vk c
    		\\&=\frac{1}{\left(\prod_{i=1}^{d}\lambda_i\right)u^{d}}e^{-\frac{(\widetilde{\vk a}u)^T\Sigma^{-1}\widetilde{\vk a}u}{2}}\int_{[0,1]^d}e^{\frac{-2u\vk c^T\vk\lambda-\vk c^T\Sigma^{-1}\vk c}{2}}d\vk c
    		\\&\sim\pk{\vk W(1)>\vk au}\int_{[0,1]^d}e^{\frac{-2u\vk c^T\vk\lambda-\vk c^T\Sigma^{-1}\vk c}{2}}d\vk c
    	\end{align*}
    	as \(u\to\IF\). Using dominated convergence theorem, one can obtain that
    	\begin{align*}
    		\int_{[0,1]^d}e^{\frac{-2u\vk c^T\vk\lambda-\vk c^T\Sigma^{-1}\vk c}{2}}d\vk c
    		=\frac{1}{u^{d}}\int_{u[0,1]^{d}}e^{\frac{-2\vk c^T\vk\lambda-\left(\vk c/\vk u\right)^T\Sigma^{-1}\left(\vk c/\vk u\right)}{2}}d\vk c
    		\sim \frac{1}{u^{d}}\int_{\R_{+}^{d}}e^{-\vk c^T\vk\lambda}d\vk c
    		=\frac{1}{\left(\prod_{i=1}^{d}\lambda_i\right)u^{d}}.
    	\end{align*}
    	Therefore we obtain that
    	\[
    	\pk{\vk W(1)>\vk au+\vk\eta}\sim \frac{1}{\left(\prod_{i=1}^{d}\lambda_i\right)u^{d}}\pk{\vk W(1)>\vk au}
    	\]
    	as \(u\to\IF\) and, by \Cref{Th1}, we have
    	\[
    	\psi_1(\vk au)\sim I_{a}u^{-d}\pk{\vk W(1)>\vk au}, \ \ \ u\to\IF,
    	\]
    	where 
    	\[
    	I_{a}=\int_{\R^{d}}\pk{\exists_{t\geq 0}:\vk W(t)-\vk at>\vk x}e^{\langle\lambda,\vk x\rangle}d\vk x.
    	\]
    \end{example}
    \section{Proofs}
    \subsection{Proof of Theorem \ref{Th1}}
    In this subsection, we formulate several auxiliary statements that will be used to prove \Cref{Th1}. All assumptions from \Cref{Th1} are valid for \Cref{Kor} and Lemmas \ref{m(u,Lambda,d)} and \ref{M(u,Lambda,d)}. As in \cite{Pdmc}, we present the upper bound for $\psi_T(\vk u)$. 
    \begin{prop} \label{Kor}
		If {$\vk u\in \R^d\setminus (-\infty,0]^d$}, then
		\BQN
		\label{fr}			\psi_T(\vk u)\leq C_{\vk K_2,T}\pk{\vk W(T)>\vk u+\vk\eta T}
		\EQN 
		for some positive constant $C_{\vk K_2,T}<\IF$.
	\end{prop}
	Let in the following \(\delta(u,\Lambda)=1-\Lambda/u^2,\) for $\Lambda>0$ and denote 
	$$m(u,\Lambda)=\pk{\exists_{t\in [0,\delta(\Lambda,u)]}:\vk W(t)-\vk\eta t>\vk au} \text{ and } M(u,\Lambda)=\pk{\exists_{ t\in [\delta(\Lambda,u),1]}:\vk W(t)-\vk\eta t>\vk au}.$$
	The following upper bound holds for $m(u,\Lambda)$.
	\begin{lem}
		For any $\Lambda>0$, $\vk a\in\R^d\setminus(-\IF,0]^d$ and a sufficiently large $u$
		\begin{align*}
			{m(u,\Lambda)} \leq C_{\vk K_2,1}e^{-\frac{\Lambda}{H}}
			\pk{\vk W(1)\ge \vk au+\vk\eta},
		\end{align*}
		where $C_{\vk K_2,1}$ is the constant from \Cref{Kor} and $H$ is a positive constant.
		\label{m(u,Lambda,d)}
	\end{lem}
For $M(u, \Lambda)$, the asymptotics as $u\to\IF$ look as follows.

	\BEL
	If $\Lambda>0$, then  
	\BQN
	M(u,\Lambda)&\sim & \left(\prod_{i\in I}\lambda_{i}\right)I_a(\Lambda)\pk{\vk W(1)>\vk au+\vk\eta}, \quad u\to \IF,
	\EQN
	with \(\vk\lambda=\Sigma^{-1}\widetilde{\vk a}\) and
	$$I_a(\Lambda)=\int_{\R^{|I|}}
	\pk{\exists_{t\in[0,\Lambda]}:\begin{array}{ccc}
			\vk W_{I}(t) - \vk a_{I}t > \vk x_{I}
	\end{array}} e^{\langle\vk \lambda_I,\vk x_I\rangle}\, \td\vk x_{I}\in (0,\IF).$$
	\label{M(u,Lambda,d)}
	\EEL
    Since \(M(u,\Lambda)\leq \psi_1(\vk au)\leq M(u,\Lambda)+m(u,\Lambda)\), using Lemmas \ref{m(u,Lambda,d)} and \ref{M(u,Lambda,d)}, we obtain that
	$$\lim_{\Lambda\to\IF}\lim_{u\to\IF}\frac{m(u,\Lambda)}{\psi_1(\vk au)}=0$$
	and
	$$\lim_{\Lambda\to\IF}\lim_{u\to\IF}\frac{M(u,\Lambda)}{\psi_1(\vk au)}=1.$$
	By the monotone convergence theorem for integrals and \Cref{Kor}, one can see that
	$$\lim_{\Lambda\to\IF}I_{a}(\Lambda)=\int_{\R^{|I|}}\pk{\exists_{t\geq 0}:\vk W_I(t)-\vk a_It>\vk x_I}e^{\langle\vk \lambda_I,\vk x_I\rangle}d\vk x_I<\IF.$$
	Therefore
	$$\psi_{1}(\vk au)\sim \left(\prod_{i\in I}\lambda_i\right)\left(\int_{\R^{|I|}}\pk{\exists_{t\geq 0}:\vk W_I(t)-\vk a_It>\vk x_I}e^{\langle\lambda_I,\vk x_I\rangle}d\vk x_I\right)\pk{\vk W(1)>\vk au+\vk\eta},\quad u\to\IF\qed .$$
	
	In the rest part of this section, we discuss the proofs of \Cref{Kor} and Lemmas \ref{m(u,Lambda,d)} and \ref{M(u,Lambda,d)}. The first subsection contains the proof of \Cref{Kor}.
	\subsection{Upper bound for \(\psi_T(\vk u)\)}
    First, let us fix the value $\vk\eta=\vk c$ and calculate
	$$\pk{\exists_{t\in [0,T]}: \vk W(t)- \vk\eta t> \vk u|\vk\eta=\vk c}.$$
	As it was shown in \cite[Thm $1.1$]{Pdmc}
	$$\pk{\exists_{t\in [0,T]}: \vk W(t)- \vk\eta t> \vk u|\vk\eta=\vk c}\leq\frac{\pk{\vk W(T)- \vk\eta T > \vk u|\vk\eta=\vk c}}{\pk{\vk W(T)> \max(\vk cT,\vk 0)|\vk\eta=\vk c}}.$$
	Our goal is now to estimate from above $$\frac{1}{\pk{\vk W(T)> \max(\vk cT,\vk 0)|\vk\eta=\vk c}}=\frac{1}{\pk{\vk W(T)> \max(\vk cT,\vk 0)}}$$ 
	or, equivalently, to find a lower bound for 
	$$\pk{\vk W(T)> \max(\vk cT,0)},$$ 
	which does not depend on the choice of values $\vk c.$ Observe that 
	\begin{align*}
		\pk{\vk W(T)> \max(\vk cT,\vk 0)}&\geq \pk{\vk W(T)> \vk K_2T}=:\frac{1}{C_{\vk K_2,T}}>0.
	\end{align*}
	Thus, we obtain that 
	$$\pk{\exists_{t\in [0,T]}:\vk W(t)- \vk\eta t> \vk u|\vk\eta=\vk c}\leq C_{\vk K_2,T}\pk{\vk W(T)-\vk\eta T > \vk u|\vk\eta=\vk c}.$$
	Since the constant $C_{K,T}$ does not depend on the choice of $\vk c=(c_1,\dots,c_d)$, from the law of total probability, it follows that
	$$\pk{\exists_{t\in [0,T]}: \vk W(t)- \vk\eta t> \vk u}\leq C_{\vk K_2,T}\pk{\vk W(T)- \vk\eta T > \vk u}.$$
	\qed
	\subsection{Proof of Lemma \ref{m(u,Lambda,d)}}
	The idea of the proof is very similar to the proof of \Cref{Kor}. Let us estimate first 
	\[
	\pk{\exists_{t\in [0,\delta(u,\Lambda)]}: \vk W(t) - \vk \eta t > \vk au|\vk\eta = \vk c}.
	\]
	By \cite[Lem \(4.6\)]{Pdmc}, one can see that
	\[
	\pk{\exists_{ t \in [0,\delta(u,\Lambda)]}:\vk W(t) - \vk \eta t > \vk au|\vk\eta = \vk c}
	\leq e^{-\frac{\Lambda}{H}}\frac{\pk{\vk W(1)- \vk\eta  > \vk au|\vk\eta = \vk c}}
	                                {\pk{\vk W(1) > \max(\vk\eta,\vk 0)|\vk\eta = \vk c}}
	\]
	for some \(H>0\) and does not depend on the choice of \(\vk c\) (see proof of \cite[Lem \(4.6\)]{Pdmc}). 
	
	As it was already shown in the proof of \Cref{Kor}
	\[
	\frac{1}{\pk{\vk W(1) > \max(\vk\eta,\vk 0)|\vk\eta = \vk c}}\leq C_{\vk K_2, 1}.
	\]
	Hence we conclude that
	\[
	\pk{\exists_{t\in [0,\delta(u,\Lambda)]}: \vk W(t) - \vk \eta t > \vk au|\vk\eta = \vk c}\leq C_{\vk K_2, 1}e^{-\frac{\Lambda}{H}}\pk{\vk W(1)- \vk\eta  > \vk au|\vk\eta = \vk c}.
	\]
	After applying the law of total probability for both hand sides we obtain the desired inequality.
	\qed
	
	We next formulate an auxiliary statement, which will play an important role in the proof of Lemma \ref{M(u,Lambda,d)}.
	
	\subsection{The uniformity.} 
    Hereafter $\phi_{\Sigma}(\vk x)$ denotes the density function of \(\vk W(1)\).
	
	In \cite[Lem \(5.1\)]{Simruin}, it was proven that
	$$\pk{\vk W(1)>\vk au+\vk c}\sim \frac{u^{-|I|}\phi_{\Sigma}(\widetilde{\vk a}u+\vk c)}{\prod_{i\in I}\lambda_i}\int_{\R^{|J|}}\mathbb{I}(\vk x_{U}\leq \vk 0_{U})e^{\langle\widetilde{\vk c}_J, \vk x_J\rangle-\frac{1}{2}\vk x_{J}^{T}(\Sigma^{-1})_{JJ}\vk x_J}\td \vk x_J,$$
	where $\vk a, \Sigma, \vk c$ are fixed, $\vk a\in\R^d\setminus(-\IF,0]^d$, $J=I^c$, $U=\{i\in J:\ \widetilde{a}_i=a_i\}$ and \(\widetilde{\vk c}=\Sigma^{-1}\vk c\). 
	
	Is this asymptotic equivalence uniform in $\vk c$? For $\vk c$ taken from a bounded set, it turns out to be true.
	\BEL
	Assume $\vk K_1\leq \vk c\leq \vk K_2$ for $\vk K_1,\vk K_2\in\R^d$ and $\ve>0$, where $<\cdot, \cdot>$ is the scalar product. If $\widetilde{\vk c}=\Sigma^{-1}\vk c$, then there exists $M=M(\ve)$ such that for $u>M$
	\label{unsim}
	\BQN
	1-\ve\leq\frac{\pk{\vk W(1)>\vk au+\vk c}\prod_{i\in I}\lambda_i}{\displaystyle{u^{-|I|}\phi_{\Sigma}(\widetilde{\vk a}u+\vk c)\int_{\R^{|J|}}\mathbb{I}(\vk x_{U}\leq \vk 0_{U})e^{\langle\widetilde{\vk c}_J, \vk x_J\rangle-\frac{1}{2}\vk x_{J}^{T}(\Sigma^{-1})_{JJ}\vk x_J}}\td \vk x_J}\leq 1+\ve,
	\EQN
	for all $\vk c$. 
	\EEL
	\begin{remark}
		In \Cref{unsim}, $M$ might however depend on $\vk a$ and $\Sigma$.
	\end{remark}
    \begin{proof}
    	The proof is based on the proof of \cite[Lem \(5.1\)]{Simruin}. For \(u>0\), introduce the \(d\)-dimensional vector \(\widetilde{\vk u}\) such that
    	\begin{align*}
    		\widetilde{\vk u}_i=\begin{cases}
    			u, \ & \text{ if } \ i\in I;
    			\\1, \ & \text{ otherwise}
    		\end{cases}
    	\end{align*}
    	for all \(i\in\{1,...,d\}\). As in \cite{Simruin}
    	\begin{align*}
    		\pk{\vk W(1)>\vk au+\vk c}&=u^{-|I|}\phi_{\Sigma}(\widetilde{\vk a}u+\vk c)\int_{\R^{d}}\mathbb{I}(\vk x<u \widetilde{\vk u}(\widetilde{\vk a}-\vk a))\theta_{u}(\vk x, \vk c)\td \vk x\\
    		&=u^{-|I|}\phi_{\Sigma}(\widetilde{\vk a}u+\vk c)\int_{\R^{d}}\mathbb{I}(\vk x_{I\cup U}<\vk 0_{I\cup U})\mathbb{I}(\vk x_{J\setminus U}<u (\widetilde{\vk a}_{J\setminus U}-\vk a_{J\setminus U}))\theta_{u}(\vk x, \vk c)\td \vk x,
    	\end{align*}
    	where 
    	$$\theta_{u}(\vk x, \vk c)=\frac{\phi_{\Sigma}(u\widetilde{\vk a}+\vk c-\frac{\vk x}{\widetilde{\vk u}})}{\phi_{\Sigma}(u\widetilde{\vk a}+\vk c)}.$$
    	The condition $\vk x_{I\cup U}<\vk 0_{I\cup U}$ implies, in particular, that we can assume $\vk x_I$ to have only negative components. Observe that
    	\begin{align*}
    		\theta_{u}(\vk x, \vk c)&=\exp\left(u\widetilde{\vk a}^T\Sigma^{-1}\left(\frac{\vk x}{\widetilde{\vk u}}\right)+\vk c^T\Sigma^{-1}\left(\frac{\vk x}{\widetilde{\vk u}}\right)-\frac{1}{2}\left(\frac{\vk x}{\widetilde{\vk u}}\right)^T\Sigma^{-1}(\frac{\vk x}{\widetilde{\vk u}})\right)\\
    		&=\exp\left(\langle\vk\lambda_I+\frac{\widetilde{\vk c}_I}{u},\vk x_I\rangle+\langle\widetilde{\vk c}_J,\vk x_J\rangle-\frac{\vk x^T_I(\Sigma^{-1})_{II}\vk x_I}{2u^2} -\frac{\vk x^T_I(\Sigma^{-1})_{IJ}\vk x_J}{u}-\frac{\vk x^T_J(\Sigma^{-1})_{JJ}\vk x_J}{2}\right).
    	\end{align*}
    	Note that since $\vk c$ is bounded from above and below, there exist constants $\widetilde{\vk c}_{I,1}$ and $\widetilde{\vk c}_{I,2}$ that do not depend on $\vk c$ such that
    	$$\widetilde{\vk c}_{I,1}\leq\widetilde{\vk c}_{I}\leq\widetilde{\vk c}_{I,2}.$$
    	If we consider $\theta_{u}(\vk x, \vk c)$ as a function $F(\widetilde{\vk c}_{I})$ of $\widetilde{\vk c}_{I}$ (which may also depend on $\vk x$, $u$ and $\widetilde{\vk c}_{J}$), then one can easily see that
    	$$F(\widetilde{\vk c}_{I,2})\leq F(\widetilde{\vk c}_{I})\leq F(\widetilde{\vk c}_{I,1}),$$
    	since we assumed that $\vk x_I$ consists of only non-positive components.
    	Hence, if the statement can be proven for the cases $\widetilde{\vk c}_{I}=\widetilde{\vk c}_{I,j}$ for $j=1,2$, the statement will be true for all possible values of $\widetilde{\vk c}_{I}$. Therefore, without loss of generality, we can assume below that $\widetilde{\vk c}_{I}$ is fixed.
    	Thus, we have that
    	$$\theta_{u}(\vk x, \vk c)=\theta(\vk x, \vk c)d_{u}(\vk x),$$
    	where $$\theta(\vk x, \vk c)=\exp\left(\langle\vk\lambda_I,\vk x_I\rangle+\langle\widetilde{\vk c}_J,\vk x_J\rangle-\frac{1}{2}\vk x^T_J(\Sigma^{-1})_{JJ}\vk x_J\right),$$ and $d_{u}(\vk x)$ is independent of $\vk c$ and converges to $1$ for all $\vk x\in \R^d$.
    	
    	Since $\vk c$ belongs to a fixed bounded set, there exists a $|J|$-dimensional vector $\vk c_{\max}$ such that
    	$$(\vk c_{\max})_j\geq|\widetilde{\vk c}_j|$$
    	for all $j\in J$ and all possible $\vk c$. Thus $$\theta(\vk x, \vk c)\leq \exp\left(\langle\vk\lambda_I,\vk x_I\rangle+\langle\vk c_{\max},|\vk x_J|\rangle-\frac{1}{2}\vk x^T_J(\Sigma^{-1})_{JJ}\vk x_J\right):=\theta_{\max}(\vk x).$$
        Since $(\Sigma^{-1})_{JJ}$ is a positive definite matrix, the function $\theta_{\max}(\vk x)$ is integrable on the set 
        
        $\{\vk x\in\R^d:\vk x_I<\vk 0_I\}$. Hence, there exists a constant $M_{\max}$ such that $$\int_{\R^d}\mathbb{I}(\vk x_{I\cup U}<\vk 0_{I\cup U})\theta(\vk x,\vk c)\td \vk x\leq M_{\max}.$$
        Using similar arguments, one can check that there exists \(M_{\min}>0\) such that
        \[
        \int_{\R^d}\mathbb{I}(\vk x_{I\cup U}<\vk 0_{I\cup U})\theta(\vk x,\vk c)\td \vk x\geq M_{\min}
        \]
        for all \(\vk c\). Therefore, it is sufficient to check that there exists $\ve_1>0$ such that $$\left|\int_{\R^{d}}\mathbb{I}(\vk x_{I\cup U}<\vk 0_{I\cup U})\mathbb{I}(\vk x_{J\setminus U}<u (\widetilde{\vk a}_{J\setminus U}-\vk a_{J\setminus U}))\theta_{u}(\vk x, \vk c)\td \vk x-\int_{\R^{d}}\mathbb{I}(\vk x_{I\cup U}<\vk 0_{I\cup U})\theta(\vk x, \vk c)\td \vk x\right|<\ve_1.$$
    	By the rectangle inequality, it is sufficient to check that for large enough $u$
    	$$\left|\int_{\R^{d}}\mathbb{I}(\vk x_{I\cup U}<\vk 0_{I\cup U})\mathbb{I}(\vk x_{J\setminus U}<u (\widetilde{\vk a}_{J\setminus U}-\vk a_{J\setminus U}))(\theta_{u}(\vk x, \vk c)-\theta(\vk x, \vk c))\td \vk x\right|<\frac{\ve_1}{2},$$
    	and
    	$$\left|\int_{\R^{d}}\mathbb{I}(\vk x_{I\cup U}<\vk 0_{I\cup U})\mathbb{I}(\vk x_{J\setminus U}> u (\widetilde{\vk a}_{J\setminus U}-\vk a_{J\setminus U}))\theta(\vk x, \vk c)\td \vk x\right|<\frac{\ve_1}{2}.$$
    	Denote here and below 
    	$$E_u:=\{\vk x\in \R^d: \vk x_{I\cup U}<\vk 0_{I\cup U}, \vk x_{J\setminus U}<u (\widetilde{\vk a}_{J\setminus U}-\vk a_{J\setminus U})\}$$ and $$E:=\{\vk x\in \R^d: \vk x_{I\cup U}<\vk 0_{I\cup U}\}.$$
    	Hence, we need to check that for large enough $u$
    	$$\left|\int_{E_u}(\theta_{u}(\vk x, \vk c)-\theta(\vk x, \vk c))\td \vk x\right|<\frac{\ve_1}{2},$$
    	and
    	$$\left|\int_{E\setminus E_u}\theta(\vk x, \vk c)\td \vk x\right|<\frac{\ve_1}{2}.$$
    	The second inequality follows from
    	$$\sup_{\vk K_1\leq\vk c\leq \vk K_2}\left|\int_{E\setminus E_u}\theta(\vk x, \vk c)\td \vk x\right|\leq \left|\int_{E\setminus E_u}\theta_{\max}(\vk x)\td \vk x\right|\to 0$$
    	as $u\to\IF$, since the sequence of sets $E_u$ converges to $E$ and $\theta_{\max}(\vk x)$ is an integrable function that does not depend on the choice of $\vk c$. Notice that 
    	$$\sup_{\vk K_1\leq\vk c\leq\vk K_2}\left|\int_{E_u}(\theta_{u}(\vk x, \vk c)-\theta(\vk x, \vk c))\td \vk x\right|\leq \sup_{\vk K_1\leq\vk c\leq\vk K_2}\int_{E}\theta(\vk x, \vk c)|d_u(\vk x)-1|\td \vk x\leq \int_{E}\theta_{\max}(\vk x)|d_u(\vk x)-1|\td \vk x,$$
    	and \(\theta_{\max}(\vk x)|d_u(\vk x)-1|\to 0\), for all $\vk x\in E$, as $u\to\IF$. To apply the dominated convergence theorem, we need to check that there exists an integrable on $E$ function $d(\vk x)$ such that
    	$$\theta_{\max}(\vk x)|d_u(\vk x)-1|\leq d(\vk x).$$
    	Since \(\theta_{\max}(\vk x)|d_u(\vk x)-1|\leq \theta_{\max}(\vk x)d_u(\vk x)+\theta_{\max}(\vk x)\)
    	and $\theta_{\max}(\vk x)$ is an integrable function on $E$, it is enough to find the corresponding function $\widetilde{d}(\vk x)$ that majorizes $\theta_{\max}(\vk x)d_u(\vk x)$ for all sufficiently large \(u\).
    	
    	Since $\Sigma^{-1}$ is positive definite, there exists $\delta>0$ such that $\Sigma^{-1}-\delta I_d$ is also positive definite.
    	Therefore, for $\vk x_I<\vk 0_I$ and sufficiently large $u$,
    	\BQNY
    	\theta_{\max}(\vk x)d_u(\vk x)&=&
    	\exp\left(\vk\lambda_I^T\vk x_I+\frac{\widetilde{\vk c}_I^T\vk x_I}{u}+\vk c_{\max}^{T}|\vk x_J|-\frac{1}{2}\left(\frac{\vk x}{\widetilde{\vk u}}\right)^T\Sigma^{-1}\left(\frac{\vk x}{\widetilde{\vk u}}\right)\right)\\
    	&\leq& \exp\left(\frac{\vk \lambda_I^T\vk x_I}{2}+\vk c_{\max}^{T}|\vk x_J|-\frac{1}{2}\left(\frac{\vk x}{\widetilde{\vk u}}\right)^T(\Sigma^{-1}-\delta I_d)\left(\frac{\vk x}{\widetilde{\vk u}}\right)-\frac{\delta}{2}\left(\frac{\vk x}{\widetilde{\vk u}}\right)^{\top}\left(\frac{\vk x}{\widetilde{\vk u}}\right)\right)\\
    	&\leq& \exp\left(\frac{\vk \lambda_I^T\vk x_I}{2}+\vk c_{\max}^{T}|\vk x_J|-\frac{\delta}{2}||\vk x_J||^2\right)=:\widetilde{d}(\vk x),
    	\EQNY
    	where \(\widetilde{d}(\vk x)\) is integrable on the set $E$ and \(\vk c_{\max}^{T}|\vk x_J|=\sum_{j\in J}(c_{\max})_{j}|x_{j}|\). Hence by the dominated convergence theorem 
    	$$\sup_{\vk K_1\leq\vk c\leq\vk K_2}\int_{E}\theta(\vk x, \vk c)|d_u(\vk x)-1|\td \vk x\leq \int_{E}\theta_{\max}(\vk x)|d_u(\vk x)-1|\td \vk x\to 0,$$
    	and
    	$$\left|\int_{E_u}(\theta_{u}(\vk x, \vk c)-\theta(\vk x, \vk c))\td \vk x\right|<\frac{\ve_1}{2},$$
    	for large $u$ and all $\vk c$ such that \(\vk K_1\leq\vk c\leq\vk K_2\). 
    	
    	Therefore as it was shown above, we obtain that for all $\ve>0$
    	\BQN
    	1-\ve\leq\frac{\pk{\vk W(1)>\vk au+\vk c}}{\displaystyle{u^{-|I|}\phi_{\Sigma}(\widetilde{\vk a}u+\vk c)\int_{\R^{d}}\mathbb{I}(\vk x_{I\cup U}\leq \vk 0_{I\cup U})e^{\langle\vk\lambda_I, \vk x_I\rangle+\langle\widetilde{\vk c}_J, \vk x_J\rangle-\frac{1}{2}\vk x_{J}^{T}(\Sigma^{-1})_{JJ}\vk x_J}}\td \vk x}\leq 1+\ve
    	\EQN
    	for $u>M(\ve)$. Finally notice that
    	\begin{align*}
    		\int_{\R^{d}}&\mathbb{I}(\vk x_{I\cup U}\leq \vk 0_{I\cup U})e^{\langle\vk\lambda_I, \vk x_I\rangle+\langle\widetilde{\vk c}_J, \vk x_J\rangle-\frac{1}{2}\vk x_{J}^{T}(\Sigma^{-1})_{JJ}\vk x_J}\td \vk x\\&=\prod_{i\in I}\left(\frac{1}{\lambda_i}\right)\int_{\R^{|J|}}\mathbb{I}(\vk x_{U}\leq \vk 0_{ U})e^{\langle\widetilde{\vk c}_J, \vk x_J\rangle-\frac{1}{2}\vk x_{J}^{T}(\Sigma^{-1})_{JJ}\vk x_J}\td \vk x_J.
    	\end{align*}
    	Therefore, we obtain
    	\BQN
    	1-\ve\leq\frac{\pk{\vk W(1)>\vk au+\vk c}\prod_{i\in I}\lambda_i}{\displaystyle{u^{-|I|}\phi_{\Sigma}(\widetilde{\vk a}u+\vk c)\int_{\R^{|J|}}\mathbb{I}(\vk x_{U}\leq \vk 0_{U})e^{\langle\widetilde{\vk c}_J, \vk x_J\rangle-\frac{1}{2}\vk x_{J}^{T}(\Sigma^{-1})_{JJ}\vk x_J}}\td \vk x_J}\leq 1+\ve,
    	\EQN
    	for large $u$. Thus the proof follows.
    \end{proof}
    
	\subsection{Proof of Lemma \ref{M(u,Lambda,d)}}
	Observe that by the law of total probability
	\[
	\pk{\exists_{t \in [\delta(u,\Lambda),1]}:\vk W(t) - \vk\eta t > \vk au}
	=\mathbb{E}_{\vk\eta}\left\{\pk{\exists_{t \in [\delta(u,\Lambda),1]}: \vk W(t) - \vk\eta t > \vk au \,\middle|\, \vk\eta = \vk c}\right\},
	\]
	where $\mathbb{E}_{\vk \eta}$ denotes the mean over $\vk\eta$. Our goal is to approximate
	\[
	\pk{\exists_{t \in [\delta(u,\Lambda),1]}: \vk W(t) - \vk\eta t > \vk au \,\middle|\, \vk\eta = \vk c}
	\]
	as \(u\to\IF\), uniformly in \(\vk c\). Since \(\vk W(t)\) is independent of \(\vk\eta\), it follows that
	\[
	\pk{\exists_{t \in [\delta(u,\Lambda),1]}: \vk W(t) - \vk\eta t > \vk au \,\middle|\, \vk\eta = \vk c}=\pk{\exists_{t \in [\delta(u,\Lambda),1]}: \vk W(t) - \vk ct > \vk au}.
	\]
	To make the presentation of the proof below easier for the reader, we consider the case
	\[
	d=2,\
	\vk a=(1,a) \ \text{ and } \
	A=\begin{pmatrix}
		1 & 0\\
		\rho & \rho^*
	  \end{pmatrix},
	\]
	where \(\rho\in(-1,1)\) and \(a\leq 1\) (as in \Cref{d=2}). Recall that the covariance matrix \(\Sigma\) of \(\vk W(1)\) is equal to 
	\(\begin{pmatrix}1 & \rho
	  \\\rho & 1\end{pmatrix}.
	\)
	
	All ideas presented in the proof can be easily generalized to the case of \(d>2\). As it was discussed above, there are only two possible cases:
	\begin{align*}
		1. \ & \ \text{ If } a>\rho, \ \text{ then } \widetilde{\vk a}=(1,a) \ \text{ and } \ I=\{1,2\};
		\\2. \ & \ \text{ If } a\leq\rho, \ \text{ then } \widetilde{\vk a}=(1,\rho) \ \text{ and } \ I=\{1\}.
	\end{align*}
	
	Consider first the case \(I=\{1,2\}\) or equivalently, \(a>\rho\) and set 
	\[
	t_u=1-\frac{t}{u^2}, \ \ \vk u_x=\vk au+\vk c-\frac{\vk x}{u}.
	\]
	Let \(\phi_{\Sigma}\) be the density function of the vector \(\vk W(1)\). Then, it follows that
	\begin{align*}
		\phi_{\Sigma}(\vk u_x)
		&=\frac{1}{2\pi\sqrt{|\det(\Sigma)|}}\exp\left(-\frac{(\vk au+\vk c-\frac{\vk x}{u})^T\Sigma^{-1}(\vk au+\vk c-\frac{\vk x}{u})}{2}\right)
		\\&=\frac{1}{2\pi\sqrt{|\det(\Sigma)|}}\exp\left(-\frac{(\vk au+\vk c)^T\Sigma^{-1}(\vk au+\vk c)}{2}+\frac{\vk x^T\Sigma^{-1}(\vk au+\vk c)}{u}-\frac{\vk x^T\Sigma^{-1}\vk x}{2u^2}\right)
		\\&=\phi_{\Sigma}(u+c_1,au+c_2)e^{\lambda_1x_1+\lambda_2x_2}\psi_u(\vk x,\vk c), \ \ \psi_u(\vk x, \vk c)=\exp\left(\frac{\vk x^T\Sigma^{-1}\vk c}{u}-\frac{\vk x^T\Sigma^{-1}\vk x}{2u^2}\right).
	\end{align*}
    where \(\psi_u(\vk x, \vk c)\to 1\), as \(u\to\IF\), for all vectors \(\vk x, \vk c.\)
    
    Our next goal is to find upper and lower bounds for \(\psi_{u}(\vk x,\vk c)\) that are uniform in \(\vk c\). To this end, observe that there exist bounded functions on \(\R^2\), \(\vk h_{\min}\) and \(\vk h_{\max}\), that do not depend on \(\vk c\) such that
	\BQN
    \label{<}\psi_u(\vk x, \vk h_{\min}(\vk x))\leq\psi_u(\vk x, \vk c)\leq\psi_u(\vk x, \vk h_{\max}(\vk x)),
	\EQN
	for all \(\vk x\in\R^2\) and all \(\vk c\in [K_{11},K_{21}]\times[K_{12},K_{22}]\). Indeed, there exist constants \(L_1,L_2,R_1,R_2\in\R\) such that 
	\[
	\Sigma^{-1}\vk c\in [L_1,R_1]\times[L_2,R_2] \ \ \ \ \ \ \ \forall \vk c\in[K_{11},K_{21}]\times[K_{12},K_{22}]
	\]
	and hence we can define \(\vk h_{\min}, \vk h_{\max}\) as follows:
	\[
	\vk h_{\min}(x_1,x_2)=\Sigma \vk q_{\min}(x_1,x_2) \ \text{ and } \ \vk h_{\max}(x_1,x_2)=\Sigma \vk q_{\max}(x_1,x_2),
	\]
	where
	\begin{align*}
	\vk q_{\min}(x_1,x_2)=\begin{cases}
    (L_1, L_2), \ & \ \text{ if } x_1>0, x_2>0;
	\\ (R_1, L_2), \ & \ \text{ if } x_1\leq 0, x_2>0;
	\\ (L_1, R_2), \ & \ \text{ if } x_1>0, x_2\leq 0;
	\\ (R_1, R_2), \ & \ \text{ if } x_1\leq 0, x_2\leq 0;
	\end{cases}
	\ \text{ and } \
    \vk q_{\max}(x_1,x_2)=\begin{cases}
    	(R_1, R_2), \ & \ \text{ if } x_1>0, x_2>0;
    	\\ (L_1, R_2), \ & \ \text{ if } x_1\leq 0, x_2>0;
    	\\ (R_1, L_2), \ & \ \text{ if } x_1>0, x_2\leq 0;
    	\\ (L_1, L_2), \ & \ \text{ if } x_1\leq 0, x_2\leq 0.
    \end{cases}
	\end{align*}
	By definition we have that \(\Sigma^{-1}\vk h_{\min}=\vk q_{\min}\) and \(\Sigma^{-1}\vk h_{\max}=\vk q_{\max}\) and, therefore, that \(\vk h_{\min}, \vk h_{\max}\) indeed satisfy inequalities (\ref{<}). By the law of total probability, one can conclude that
	\begin{align*}
		\pk{\exists_{t \in [\delta(u,\Lambda),1]}:\vk W(t) - \vk ct > \vk au}=\frac{1}{u^2}\int_{\R^2}\pk{\exists_{t \in [0,\Lambda]}: \vk W(t_u)- \vk c t_u  > \vk au \,\middle|\,
			    \vk W(1)=\vk u_x
				}\phi_{\Sigma}(\vk u_x)d\vk x
		\\=\frac{1}{u^2}\phi_{\Sigma}(\vk au+\vk c)\int_{\R^2}\pk{\exists_{t \in [0,\Lambda]}: \vk W(t_u)- \vk c t_u  > \vk au \,\middle|\,
			\vk W(1)=\vk u_x
		}e^{\lambda_1x_1+\lambda_2x_2}\psi_{u}(\vk x,\vk c)d\vk x.
	\end{align*}
	
	Let us prove that 
	\[
	I_{a,u}(\vk c, \Lambda)=\int_{\R^2}\pk{\exists_{t \in [0,\Lambda]}: \vk W(t_u)- \vk c t_u  > \vk au \,\middle|\,
		\vk W(1)=\vk u_x
	}e^{\lambda_1x_1+\lambda_2x_2}\psi_{u}(\vk x,\vk c)d\vk x\to I_{a}(\Lambda)
	\]
	uniformly in \(\vk c\). To start with, observe that
	\begin{align*}
		L_{a, u}(\vk c, \Lambda)&=\int_{\R^2}\pk{\exists_{t \in [0,\Lambda]}: \vk W(t_u)- \vk c t_u  > \vk au \,\middle|\,
			\vk W(1)=\vk u_x
		}e^{\lambda_1x_1+\lambda_2x_2}\psi_{u}(\vk x,\vk h_{\min}(\vk x))d\vk x
		\\&\leq\int_{\R^2}\pk{\exists_{t \in [0,\Lambda]}: \vk W(t_u)- \vk c t_u  > \vk au \,\middle|\,
			\vk W(1)=\vk u_x
		}e^{\lambda_1x_1+\lambda_2x_2}\psi_{u}(\vk x,\vk c)d\vk x
		\\&\leq\int_{\R^2}\pk{\exists_{t \in [0,\Lambda]}: \vk W(t_u)- \vk c t_u  > \vk au \,\middle|\,
			\vk W(1)=\vk u_x
		}e^{\lambda_1x_1+\lambda_2x_2}\psi_{u}(\vk x,\vk h_{\max}(\vk x))d\vk x=R_{a, u}(\vk c, \Lambda).
	\end{align*}
	It turns out that \(L_{a, u}(\vk c, \Lambda)\) and \(R_{a, u}(\vk c, \Lambda)\) are independent of \(\vk c\), and hence if we check that
	\BQN
	\label{h_min,max}\lim_{u\to\IF}L_{a, u}(\vk c, \Lambda)=\lim_{u\to\IF}R_{a, u}(\vk c, \Lambda)=I_{a}(\Lambda),
	\EQN
	
	then the desired uniform convergence \(I_{a,u}(\vk c,\Lambda)\to I_{a}(\Lambda)\) follows too. 
	
	Thus, it is enough to check (\ref{h_min,max}). We will prove that \(\lim_{u\to\IF}R_{a, u}(\vk c, \Lambda)=I_{a}(\Lambda)\), other equality can be proven in a similar way. Denote by \(\vk D(t)=\vk W(t)-t\vk W(1)\). One can check that 
    \[
    cov(\vk D(t), \vk W(1))=cov(\vk W(t), \vk W(1))-tcov(\vk W(1), \vk W(1))=\vk 0_{d\times d}
    \]
    for all \(t\geq 0\), with \(\vk 0_{d\times d}\) being a \(d\times d\) zero matrix. Here by \(cov(\vk w_1,\vk w_2)\), with \(d\)-dimensional random vectors \(\vk w_1,\vk w_2\), we denote the covariance matrix of \(\vk w_1,\vk w_2\).
	Since \(\vk D\) is a Gaussian process, we conclude that the process \(\vk D\) is independent of the vector \(\vk W(1)\). Therefore, the event
	\[
    \vk W(t_u)- \vk ct_u  > \vk au \ \Big | \ \vk W(1)=\vk au+\vk c-\frac{\vk x}{u}.
	\]
	is equivalent to
	\begin{align*}
		\vk D(t_u) + t_u\vk u_x- \vk ct_u  > \vk au
		&\iff u\left(\vk D(t_u) + \left(1-\frac{t}{u^2}\right)\left(\vk au-\frac{\vk x}{u}\right) - \vk au\right)>\vk 0
		\\&\iff \vk X_u(t)=u\vk D(t_u)-\vk at-\vk x+\frac{t\vk x}{u^2}>\vk 0.
	\end{align*}
	According to the proof of \cite[Lem \(3.2\)]{mi:18} the process \(\vk X_u(t)\) weakly converges to \(\vk W(t)-\vk at-\vk x\), as \(u\to\IF\).
	
	Moreover, \(\vk X_u(t)\) does not depend on \(\vk c\), hence \(L_{a, u}(\vk c, \Lambda)\) and \(R_{a, u}(\vk c, \Lambda)\) do not depend on \(\vk c\) too and the weak convergence is uniform in \(\vk c\).
	
	Since the supremum functional is continuous on \(C\left([0,\Lambda]\right)\), using continuous mapping theorem, one can conclude that
	\begin{align*}
		\pk{\sup_{t\in[0,\Lambda]}\min\left(uD_1(t_u)-t-x_1+\frac{tx_1}{u^2},uD_2(t_u)-at-x_2+\frac{tx_2}{u^2}\right)>0}
		\\ \to\pk{\sup_{t\in[0,\Lambda]}\min\left(W_1(t)-t-x_1, W_2(t)-at-x_2\right)>0}, \ \ \ u\to\IF
	\end{align*}
	for almost all \(\vk x\in\R^2\) uniformly in \(\vk c\). Therefore
    \[
    \pk{\exists_{t \in [0,\Lambda]}: \vk X_u(t)>\vk 0}e^{\lambda_1x_1+\lambda_2x_2}\psi_{u}(\vk x,\vk h_{\max}(\vk x))\to\pk{\exists_{t \in [0,\Lambda]}:\vk W(t)-\vk at-\vk x>\vk 0}e^{\lambda_1x_1+\lambda_2x_2}
    \]
    almost everywhere, as \(u\to\IF\). To apply the dominated convergence theorem we need to check that there exists an integrable function \(F(\vk x)\) on \(\R^2\) such that
    \BQN
    \label{F}\pk{\exists_{t \in [0,\Lambda]}: \vk X_u(t)>\vk 0}e^{\lambda_1x_1+\lambda_2x_2}\psi_{u}(\vk x,\vk h_{\max}(\vk x))\leq F(\vk x).
    \EQN
    for all large enough \(u\). In the proof of \cite[Lem \(3.2\)]{mi:18} it was shown that there exist constants \(\overline{C}, C\) that are independent of the choice of \(\vk c\) such that
    \begin{align*}
    	\pk{\exists_{t \in [0,\Lambda]}: \vk X_u(t)>\vk 0}
    	\leq\begin{cases}
    		\overline{C}e^{-C(x_1^2+x_2^2)}, \ & \ x_1>0,x_2>0;
    		\\ \overline{C}e^{-Cx_2^2}, \ & \ x_1\leq 0,x_2>0;
    		\\ \overline{C}e^{-Cx_1^2}, \ & \ x_1>0,x_2\leq 0;
    		\\ 1, \ & \ x_1\leq 0, x_2\leq 0.
    	\end{cases}
    \end{align*}
    Define the function \(F(x_1,x_2)\) as follows:
    \begin{align*}
    	F(x_1,x_2)
    	=\begin{cases}
    		\overline{C}e^{-C(x_1^2+x_2^2)}e^{2\lambda_1x_1+\lambda_2x_2}, \ & \ x_1>0,x_2>0;
    		\\ \overline{C}e^{-Cx_2^2}e^{\frac{\lambda_1}{2}x_1+2\lambda_2x_2}, \ & \ x_1\leq 0,x_2>0;
    		\\ \overline{C}e^{-Cx_1^2}e^{2\lambda_1x_1+\frac{\lambda_2}{2}x_2}, \ & \ x_1>0,x_2\leq 0;
    		\\ e^{\frac{\lambda_1}{2}x_1+\frac{\lambda_2}{2}x_2}, \ & \ x_1\leq 0, x_2\leq 0.
    	\end{cases}
    \end{align*}
    One can check that this function is integrable and satisfy (\ref{F}).
    
    Then by the dominated convergence theorem we obtain that \(R_{a, u}(\vk c, \Lambda)\to I_{a}(\Lambda)\), as \(u\to\IF\). 
    
    By the similar reasons \(L_{a, u}(\vk c, \Lambda)\to I_{a}(\Lambda)\), as \(u\to\IF\). Hence
    \[
    I_{a,u}(\vk c, \Lambda)\to I_{a}(\Lambda),
    \]
	as \(u\to\IF\), uniformly in \(\vk c\). Since \(I=\{1,2\}\), \(J=U=\emptyset\) and therefore by \Cref{unsim} we obtain that
	\[
    \frac{1}{u^{2}}\phi_{\Sigma}(\widetilde{\vk a}u+\vk c)\sim \lambda_1\lambda_2\pk{\vk W(1)>\vk au+\vk c},
	\]
	and that the asymptotic equivalence is uniform in \(\vk c\in[K_{11},K_{21}]\times[K_{12},K_{22}]\). Hence, we have that if \(I=\{1,2\}\), then
	\[
	\pk{\exists_{t \in [\delta(u,\Lambda),1]}: \vk W(t) - \vk\eta t > \vk au \,\middle|\, \vk\eta = \vk c}\sim \lambda_1\lambda_2\pk{\vk W(1)>\vk au+\vk c}I_{a}(\Lambda),
	\]
	uniformly in \(\vk c\). Thus, by the law of total probability, we obtain that
	\[
	M(u,\Lambda)=\mathbb{E}_{\vk\eta}\left\{\pk{\exists_{t \in [\delta(u,\Lambda),1]}: \vk W(t) - \vk\eta t > \vk au \,\middle|\, \vk\eta = \vk c}\right\}\sim\lambda_1\lambda_2\pk{\vk W(1)>\vk au+\vk \eta}I_{a}(\Lambda),
	\]
	that finishes the proof of Lemma \ref{M(u,Lambda,d)}, if \(a>\rho\).
	
	Consider now the case \(a\leq \rho\).  As above denote 
	\[
	t_u:=1-\frac{t}{u^2}, \vk u_x:=(u+c_1-x_1/u, \rho u+c_2-x_2) \ \text{ and } \ \widetilde{\vk u}:=(u,1).
	\]
	Notice that
	\begin{align*}
		\phi_{\Sigma}(\vk u_x)
		&=\frac{1}{2\pi\sqrt{|\det(\Sigma)|}}\exp\left(-\frac{(\widetilde{\vk a}u+\vk c-\vk x/\widetilde{\vk u})^T\Sigma^{-1}(\widetilde{\vk a}u+\vk c-\vk x/\widetilde{\vk u})}{2}\right)
		\\&=\frac{1}{2\pi\sqrt{|\det(\Sigma)|}}\exp\left(-\frac{(\widetilde{\vk a}u+\vk c)^T\Sigma^{-1}(\widetilde{\vk a}u+\vk c)}{2}+\left(\frac{\vk x}{\widetilde{\vk u}}\right)^T\Sigma^{-1}(\widetilde{\vk a}u+\vk c)-\frac{1}{2}\left(\frac{\vk x}{\widetilde{\vk u}}\right)^T\Sigma^{-1}\left(\frac{\vk x}{\widetilde{\vk u}}\right)\right)
		\\&=:\phi_{\Sigma}(\widetilde{\vk a}u+\vk c)\theta_u(\vk x, \vk c).
	\end{align*}
	One can check that
	\begin{align*}
		\theta_u(\vk x, \vk c)
		&=\exp\left(\left(\frac{x_1}{u}\right)(\Sigma^{-1}(\widetilde{\vk a}u+\vk c))_1+x_2(\Sigma^{-1}(\widetilde{\vk a}u+\vk c))_2-\frac{1}{2}\left(\frac{\vk x}{\widetilde{\vk u}}\right)^T\Sigma^{-1}\left(\frac{\vk x}{\widetilde{\vk u}}\right)\right)
		\\&=\exp\left(x_1+\frac{(c_1-\rho c_2)x_1}{u(1-\rho^2)}+\frac{(c_2-\rho c_1)x_2}{(1-\rho^2)}-\frac{x_1^2}{2(1-\rho^2)u^2}+\frac{\rho x_1x_2}{u(1-\rho^2)}-\frac{x_2^2}{2(1-\rho^2)}\right)
		\\&=:\exp\left(x_1-\frac{x_2^2-2(c_2-\rho c_1)x_2}{2(1-\rho^2)}\right)\psi_u(\vk x, \vk c),
	\end{align*}
	where \(\lim_{u\to\IF} \psi_u(\vk x, \vk c)=1\) for all \(\vk x, \vk c\in\R^2.\)
	Thus we have
	\[
	\phi_{\Sigma}(\vk u_x)=\phi_{\Sigma}(\widetilde{\vk a}u+\vk c)\exp\left(x_1-\frac{x_2^2-2(c_2-\rho c_1)x_2}{2(1-\rho^2)}\right)\psi_u(\vk x, \vk c).
	\]
	In addition, similarly to the case \(a>\rho\) one can define bounded on \(\R^2\) functions \(\vk h_{\min}(\vk x), \vk h_{\max}(\vk x)\) such that
    \BQN
    \label{h_min,max,2}\psi_u(\vk x, \vk h_{\min}(\vk x))\leq\psi_u(\vk x, \vk c)\leq\psi_u(\vk x, \vk h_{\max}(\vk x)).
    \EQN
    
    Again, without loss of generality we can replace \(\psi_u(\vk x, \vk c)\) with \(\psi_u(\vk x, \vk h_{\max}(\vk x))\). 
    
    By the law of total probability
	\begin{align*}
		\pk{\exists_{t \in [\delta(u,\Lambda),1]}:\vk W(t) - \vk ct > u}=\frac{1}{u}\int_{\R^2}\pk{\exists_{t \in [0,\Lambda]}: \vk W(t_u)- \vk c t_u  > \vk au \,\middle|\,
			\vk W(1)=\vk u_x
		}\phi_{\Sigma}(\vk u_x)d\vk x
		\\=\frac{1}{u}\phi_{\Sigma}(\widetilde{\vk a}u+\vk c)\int_{\R^2}\pk{\exists_{t \in [0,\Lambda]}: \vk W(t_u)- \vk c t_u  > \vk au \,\middle|\,
			\vk W(1)=\vk u_x
		}e^{x_1-\frac{x_2^2-2(c_2-\rho c_1)x_2}{2(1-\rho^2)}}\psi_u(\vk x, \vk c)d\vk x.
	\end{align*}
	The event 
	\[
	\vk W(t_u)- \vk ct_u  > \vk au \ \Big | \ \vk W(1)=\widetilde{\vk a}u+\vk c-\frac{\vk x}{\widetilde{\vk u}}
	\]
	is equivalent to
	\begin{align*}
		\vk D(t_u) + t_u\vk u_x- \vk ct_u  > \widetilde{\vk a}u
		&\iff u\left(\vk D(t_u) + \left(1-\frac{t}{u^2}\right)\left(\widetilde{\vk a}u-\frac{\vk x}{\widetilde{\vk u}}\right) - \vk au\right)>\vk 0
		\\&\iff \vk X_u(t)=u\vk D(t_u)-\widetilde{\vk a}t-\frac{u\vk x}{\widetilde{\vk u}}+\frac{t\vk x}{u\widetilde{\vk u}}+(\widetilde{\vk a}-\vk a)u^2>\vk 0
		\\&\iff\begin{array}{ccc}
			uD_1(t_u)-t-x_1+\frac{tx_1}{u^2}>0; \\
			uD_2(t_u)-\rho t-ux_2+\frac{tx_2}{u}+(\rho-a)u^2>0.
		\end{array}
	\end{align*}
	The process \(\vk X_u(t)\) is independent of \(\vk c\), and one can check that \(\vk X_u(t)\) weakily converges to
	\begin{align*}
	& 1.\ (W_1(t)-t-x_1, \IF), \ & \text{ if } \ & \ \rho>a \ \text{ or } \ \rho=a \ \text{ and } \ x_2<0;
	\\& 2.\ (W_1(t)-t-x_1, W_2(t)-\rho t), \ & \text{ if } \ & \ \rho=a \ \text{ and } \ x_2=0;
	\\& 3.\ (W_1(t)-t-x_1, -\IF), \ & \text{ if } \ & \ \rho=a \ \text{ and } \ x_2>0;
	\end{align*}
	as \(u\to\IF\). Hence by continuous mapping theorem, we have that
	\begin{align*}
	    \pk{\exists_{t\in [0,\Lambda]}:\vk X_u(t)>\vk 0}\to
	     \begin{cases}
	     	\pk{\exists_{t\in [0,\Lambda]}:W_1(t)-t-x_1>0}, \ & \ \text{ if } \rho>a \ \text{ or } \ \rho=a \ \text{ and } \ x_2<0,
	     	\\ \pk{\exists_{t\in [0,\Lambda]}:\begin{array}{ccc}
	     			W_1(t)-t-x_1>0 \\
	     			W_2(t)-\rho t>0
	     		\end{array}}, \ & \ \text{ if } \ \rho=a \ \text{ and } \ x_2=0,
	     	\\ 0, \ & \ \text{ otherwise}
	     \end{cases}
	\end{align*}
	almost everywhere, as \(u\to\IF\). Denote by
	\[
	\alpha_u(a, \vk x)=\pk{\exists_{t\in [0,\Lambda]}:\vk X_u(t)>\vk 0}\psi_u(\vk x, \vk h_{\max}(\vk x)),
	\]
	and by
	\begin{align*}
		\alpha(\vk x)=\begin{cases}
			\pk{\exists_{t\in [0,\Lambda]}:W_1(t)-t-x_1>0}, \ & \ \text{ if } \ a<\rho;
			\\ \pk{\exists_{t\in [0,\Lambda]}:W_1(t)-t-x_1>0}\mathbb{I}(x_2<0), \ & \ \text{ if } \ a=\rho.
		\end{cases}
	\end{align*}
	Therefore, we obtain that for almost all \(\vk x\in\R^2\)
	\[
	\alpha_u(a, \vk x)e^{x_1-\frac{x_2^2-2(c_2-\rho c_1)x_2}{2(1-\rho^2)}}\to \alpha(\vk x)e^{x_1-\frac{x_2^2-2(c_2-\rho c_1)x_2}{2(1-\rho^2)}}.
	\]
    Since \(\vk c\) belongs to a bounded set, the described convergence is uniform in \(\vk c\).
	
	To finish the proof check that 
	\[
	\frac{\Big|\int_{\R^2}(\alpha_u(a,\vk x)-\alpha(\vk x))e^{x_1-\frac{x_2^2-2(c_2-\rho c_1)x_2}{2(1-\rho^2)}}d\vk x\Big|}{\int_{\R^2}\alpha(\vk x)e^{x_1-\frac{x_2^2-2(c_2-\rho c_1)x_2}{2(1-\rho^2)}}d\vk x}\to 0
	\]
    as \(u\to\IF\) uniformly in \(\vk c\). To this end we will prove that 
    \[
    \Big|\int_{\R^2}(\alpha_u(a,\vk x)-\alpha(\vk x))e^{x_1-\frac{x_2^2-2(c_2-\rho c_1)x_2}{2(1-\rho^2)}}d\vk x\Big|\to 0,
    \]
    uniformly in \(\vk c\), and check that there exists a constant \(C_{\min}>0\) such that 
    \[
    \int_{\R^2}\alpha(\vk x)e^{x_1-\frac{x_2^2-2(c_2-\rho c_1)x_2}{2(1-\rho^2)}}d\vk x>C_{\min},
    \]
    for all \(\vk c\in[K_{11}, K_{21}]\times[K_{12}, K_{22}]\). We start with the former property. To this end observe that there exists a bounded function \(K(x_2)\) such that
    \[
    \frac{(c_2-\rho c_1)x_2}{1-\rho^2}\leq K(x_2)x_2.
    \]
    Therefore, it is sufficient to check that 
    \[
    \Big|\int_{\R^2}(\alpha_u(a,\vk x)-\alpha(\vk x))e^{x_1-\frac{x_2^2-2(c_2-\rho c_1)x_2}{2(1-\rho^2)}}d\vk x\Big|\leq\int_{\R^2}|\alpha_u(a,\vk x)-\alpha(\vk x)|e^{x_1+K(x_2)x_2-\frac{x_2^2}{2(1-\rho^2)}}d\vk x\to 0,
    \]
    but it follows from the dominated convergence theorem. Indeed, we have shown above that
    \[
    |\alpha_u(a,\vk x)-\alpha(\vk x)|e^{x_1+K(x_2)x_2-\frac{x_2^2}{2(1-\rho^2)}}\to 0.
    \]
    Moreover, 
    \[
    |\alpha_u(a,\vk x)-\alpha(\vk x)|e^{x_1+K(x_2)x_2-\frac{x_2^2}{2(1-\rho^2)}}
    \leq \alpha_u(a,\vk x)e^{x_1+K(x_2)x_2-\frac{x_2^2}{2(1-\rho^2)}}
    +\alpha(\vk x)e^{x_1+K(x_2)x_2-\frac{x_2^2}{2(1-\rho^2)}}.
    \]
    Thus, it is enough to find an integrable on \(\R^2\) function \(F\) such that
    \[
    \alpha_u(a,\vk x)e^{x_1+K(x_2)x_2-\frac{x_2^2}{2(1-\rho^2)}}\leq F(\vk x),
    \]
    for large enough \(u\) and check that \(\alpha(\vk x)e^{x_1+K(x_2)x_2-\frac{x_2^2}{2(1-\rho^2)}}\) is an integrable function.
    
    We check the first statement. Recall that
    \begin{align*}
    	\alpha_u(a,\vk x)e^{x_1+K(x_2)x_2-\frac{x_2^2}{2(1-\rho^2)}}
    	&=\pk{\exists_{t\in [0,\Lambda]}:\vk X_u(t)>\vk 0}\psi_u(\vk x, \vk h_{\max}(\vk x))e^{x_1+K(x_2)x_2-\frac{x_2^2}{2(1-\rho^2)}}
    	\\&=\pk{\exists_{t\in [0,\Lambda]}:\vk X_u(t)>\vk 0}e^{x_1(1+\frac{(\vk h_{\max}(x_1))_1}{u})+K(x_2)x_2}e^{-\frac{1}{2}\left(\frac{\vk x}{\widetilde{\vk u}}\right)^T\Sigma^{-1}\left(\frac{\vk x}{\widetilde{\vk u}}\right)},
    \end{align*}
    where \((\vk h_{\max}(x_1))_1\) is the first coordinate of the vector function \((\vk h_{\max}(x_1))_1\). Since \(\Sigma^{-1}\) is a positive definite matrix, there exists \(\delta>0\) such that \(\Sigma^{-1}-\delta I_{d}\) is positive definite too. Hence, we have that
    \begin{align*}
    \alpha_u(a,\vk x)e^{x_1+K(x_2)x_2-\frac{x_2^2}{2(1-\rho^2)}}
    &\leq\pk{\exists_{t\in [0,\Lambda]}:\vk X_u(t)>\vk 0}e^{x_1(1+\frac{(\vk h_{\max}(x_1))_1}{u})+K(x_2)x_2}e^{-\frac{\delta x_1^2}{u^2}-\delta x_2^2}
    \\&\leq\pk{\exists_{t\in [0,\Lambda]}:\vk X_u(t)>\vk 0}e^{x_1(1+\frac{(\vk h_{\max}(x_1))_1}{u})+K(x_2)x_2-\delta x_2^2}.
    \end{align*}
    If \(x_1>0\), then using Borel's inequality \cite[Lem \(4.5\)]{Vector} one can check that there exist positive constants \(\overline{C}\) and \(C\) such that for large enough \(u\)
    \[
    \pk{\exists_{t\in [0,\Lambda]}:\vk X_u(t)>\vk 0}\leq \pk{\exists_{t\in [0,\Lambda]}:uD_1(t_u)-t-x_1+\frac{tx_1}{u^2}>0}\leq\overline{C}e^{-Cx_1^2}.
    \]
    Hence, one can define the function \(F\) as follows:
    \begin{align*}
    	F(x_1,x_2)=\begin{cases}
    		\overline{C}e^{-Cx_1^2+2x_1+K(x_2)x_2-\delta x_2^2}, \ & \text{ if } x_1>0;
    		\\e^{\frac{x_1}{2}+K(x_2)x_2-\delta x_2^2}, \ & \text{ otherwise.}
    	\end{cases}
    \end{align*}
    It is easy to see that \(F\) is integrable on \(\R^2\) and that
    \[
    \alpha_u(a,\vk x)e^{x_1+K(x_2)x_2-\frac{x_2^2}{2(1-\rho^2)}}\leq F(\vk x).
    \]
    The function \(\alpha(\vk x)e^{x_1+K(x_2)x_2-\frac{x_2^2}{2(1-\rho^2)}}\) is integrable by the similar reasons, using Borel's inequality, so we omit the proof here. Thus 
    \[
    \Big|\int_{\R^2}(\alpha_u(a,\vk x)-\alpha(\vk x))e^{x_1-\frac{x_2^2-2(c_2-\rho c_1)x_2}{2(1-\rho^2)}}d\vk x\Big|\to 0,
    \]
    uniformly in \(\vk c\in[K_{11}, K_{21}]\times[K_{12}, K_{22}]\). In addition, there exists a constant \(\overline{K}\) such that
    \[
    \frac{(c_2-\rho c_1)}{1-\rho^2}<\overline{K}.
    \]
    Therefore, we obtain that, since \(\alpha(\vk x)=1\) for \(x_1<0,x_2<0\), then
    \begin{align*}
    	\int_{\R^2}\alpha(\vk x)e^{x_1-\frac{x_2^2-2(c_2-\rho c_1)x_2}{2(1-\rho^2)}}d\vk x
    	&\geq\int_{x_1<0,x_2<0}e^{x_1-\frac{x_2^2-2(c_2-\rho c_1)x_2}{2(1-\rho^2)}}d\vk x
    	\\&\geq\int_{x_1<0,x_2<0}e^{x_1+\overline{K}x_2-\frac{x_2^2}{2(1-\rho^2)}}d\vk x=:C_{\min}>0.
    \end{align*}
    Hence, we can conclude that indeed
    \[
    \frac{\Big|\int_{\R^2}(\alpha_u(a,\vk x)-\alpha(\vk x))e^{x_1-\frac{x_2^2-2(c_2-\rho c_1)x_2}{2(1-\rho^2)}}d\vk x\Big|}{\int_{\R^2}\alpha(\vk x)e^{x_1-\frac{x_2^2-2(c_2-\rho c_1)x_2}{2(1-\rho^2)}}d\vk x}\to 0,
    \]
    uniformly in \(\vk c\) or, equivalently,
    \[
    \int_{\R^2}\alpha_u(a,\vk x)e^{x_1-\frac{x_2^2-2(c_2-\rho c_1)x_2}{2(1-\rho^2)}}d\vk x
    \sim\int_{\R^2}\alpha(\vk x)e^{x_1-\frac{x_2^2-2(c_2-\rho c_1)x_2}{2(1-\rho^2)}}d\vk x,
    \]
    uniformly in \(\vk c\), as \(u\to\IF\). Since 
    \begin{align*}
    	\int_{\R^2}\alpha(\vk x)e^{x_1-\frac{x_2^2-2(c_2-\rho c_1)x_2}{2(1-\rho^2)}}d\vk x
    	=\int_{\R}\pk{\exists_{t\in [0,\Lambda]}:W_1(t)-t-x_1>0}e^{x_1}dx_1\int_{\vk x_U<\vk 0_U}e^{-\frac{x_2^2-2(c_2-\rho c_1)x_2}{2(1-\rho^2)}}dx_2,
    \end{align*}
    by  \Cref{unsim} we can conclude that
    \[
    \pk{\exists_{t \in [\delta(u,\Lambda),1]}: \vk W(t) - \vk\eta t > \vk au \,\middle|\, \vk\eta = \vk c}\sim\left(\prod_{i\in I}\lambda_{i}\right)I_a(\Lambda)\pk{\vk W(1)>\vk au+\vk c}
    \]
    uniformly in \(\vk c\). Here we used that \(U=\emptyset\) if \(a<\rho\) and \(U=\{2\}\) if \(a=\rho\).
    
    Hence by the law of total probability, we obtain that
    \[
    M(u,\Lambda)\sim\left(\prod_{i\in I}\lambda_{i}\right)I_a(\Lambda)\pk{\vk W(1)>\vk au+\vk\eta}, \quad u\to \IF.
    \]
    That finishes the proof of Lemma \ref{M(u,Lambda,d)}. 
    \qed
    \section*{Acknowledgments}
    Partial  supported by SNSF Grant 200021-196888 is kindly acknowledged.
    \def\polhk#1{\setbox0=\hbox{#1}{\ooalign{\hidewidth
    			\lower1.5ex\hbox{`}\hidewidth\crcr\unhbox0}}}

\end{document}